Dmytro Taranovsky
March 10, 2012
Modified: December 14, 2016

# Reflective Cardinals

**Alternate Title:** *Higher Order Set Theory with Reflective Cardinals*

**Abstract:** We introduce and axiomatize the notion of a reflective cardinal, use it to give semantics to higher order set theory, and explore connections between the notion of reflective cardinals and large cardinal axioms.

**Contents:**

# 1 Introduction

Despite vast advances in set theory and mathematics in general, the language of set theory, which is the de facto formal language of mathematics, has remained the same since the beginning of modern set theory and first order logic. That formal language has served us well, but it has a major limitation: It does not adequately deal with properties for which there is no set of all objects satisfying the property. In particular, it does not allow quantification over such predicates, nor creation of collections of such predicates. Our purpose is to address this deficiency, and to extend the expressive power of set theory to the furthest heights.

A natural way to extend a language is to augment it with higher types, and it is not technically difficult to come up with a higher order set theory. Examples include Kelley-Morse set theory and Ackermann set theory. The problem is rather a semantic one. If a set is an arbitrary collection of objects, what sense do we make of "collections" that are not sets?

We avoid this problem by not including "collections" that are not sets. Our addition to the language of set theory is instead a predicate on sets: $R(\kappa)$ iff $\kappa$ is a reflective cardinal. We provide semantics for $R$ both by "defining" $R$ in the human language, and by creating formal axiomatic systems that include $R$. We then interpret higher order set theory in terms of ordinary set theory augmented with $R$.

In the second section, we define $R$, and interpret higher order set theory using $R$. In the third section, we axiomatize $R$ and explain the axioms. In the fourth section, we demonstrate a close connection between $R$ and large cardinal axioms. In the fifth section, we iterate the notion of a reflective cardinal, axiomatize the resulting extensions, and again demonstrate a connection with large cardinal axioms. In the sixth section, we use reflective cardinals to argue for existence of supercompact cardinals. In the seventh and eighth sections, we introduce notions of reflective sequences of ordinals — notions that go beyond the expressive power of reflective cardinals.

By its nature, the paper includes both formal results and philosophical arguments. The formal results are valid theorems independently of whether the reader agrees with us that the notion of a reflective cardinal is a well-defined one.

For background about large cardinals, the reader is referred to a standard text on set theory such as (Jech 2006), or to a text on large cardinals, such as (Kanamori 2008). For background about subtle and ineffable cardinals, see (Friedman 2001). Reflective cardinals (under the name of reflective ordinals) were originally introduced in (Taranovsky 2005). The exposition of reflective cardinals there can be seen as a precursor to this paper.

## 2 Definition and Philosophy

**Definition:** $\kappa$ is a *reflective cardinal*, denoted by $R(\kappa)$, iff $(V, \in, \kappa)$ has the same theory with parameters in $V_\kappa$ as $(V, \in, \lambda)$ for every cardinal $\lambda > \kappa$ with sufficiently strong reflection properties.
**Semantics of higher order set theory:** A statement in higher order set theory (such as a second order statement about $V$) is true iff it is true about $V_\kappa$ for a reflective $\kappa$. If $\varphi(s_1, \ldots, s_n)$ is a formula in higher order set theory and $s_1, \ldots, s_n$ are sets, then $\varphi(s_1, \ldots, s_n)$ is true iff it is true about $V_\kappa$ for a reflective $\kappa$ with $(s_1, \ldots, s_n) \in V_\kappa$.

The definition is not formalizable in the language of set theory: A proper extension to the language of set theory cannot be made in the language of set theory. Instead, just like the notion of set itself, the definition is made in the human language. We then use general reasoning to come up with axioms about the notion, and the axioms are in turn used to prove formal theorems. In turn, the theorems can provide insights that may lead to new axioms.

Informally, a reflection property states that the ordinal is in a certain way similar to Ord (the class of all ordinals). An example is: $V_\kappa$ is a $\Sigma_n$ elementary substructure of $V$. For any particular finite $n$, this is formalizable in the language of set theory, and existence of such $\kappa$ provable in ZFC.

Large cardinal properties are canonical examples of reflection properties. Large cardinal axioms are often motivated by a consideration about what properties should Ord satisfy (such as being regular), and they postulate existence of cardinals with some of these properties. Reflective cardinals satisfy all genuine large cardinal properties that are expressible in the ordinary language of set theory and are fully realized in $V$.

The definition of reflective cardinals relies on some notion of a sufficiently strong reflection property. Evidence for existence of such properties comes from inner model theory. Let $M$ be a canonical inner model (such as $L$), and let $M^\#$ be the sharp of $M$. $M^\#$ generates a proper class of indiscernibles for $M$. Inside $M$, these indiscernibles satisfy sufficiently strong reflection properties, so if $M^\#$ exists, the notion of reflective cardinals for $M$ is not problematic (and $R^M$ is definable in $V$). For $V$, however, the analogue of $V^\#$ is not definable in the language of set theory. Instead, we are postulating existence of sufficiently strong reflection properties for the definition of reflective cardinals to make sense.

By the reflection principle, in so far as a statement $\varphi$ in higher order set theory is meaningful and true, it must be true about $V_\kappa$ for some cardinal $\kappa$. If it were otherwise, then the set-theoretical universe would appear incomplete since there is new structure at the level of Ord that is not there earlier in a similar form. Since in a sufficient sense, reflective cardinals are those that are most like Ord, $\varphi$ must also be true about reflective cardinals. Now, we turn this around and use reflective cardinals to *define* semantics of higher order theory: If $\varphi(s_1, \ldots, s_n)$ is a formula in higher order set theory and $s_1, \ldots, s_n$ are sets, then $\varphi(s_1, \ldots, s_n)$ is true iff it is true about $V_\kappa$ for a reflective $\kappa$ with $(s_1, \ldots, s_n) \in V_\kappa$. This definition is correct in so far as higher order set theory is meaningful, and since the definition gives unambiguous truth value semantics, there is no sense in considering higher order set theory to be meaningless. For a reflective $\kappa$, sets correspond to sets in $V_\kappa$, and proper classes and beyond correspond to sets outside of $V_\kappa$.

To encapsulate higher order set theory with arbitrary set parameters, the definition of reflective cardinals includes two things:
* We use a predicate for reflective cardinals instead of a constant for the least reflective cardinal.
* We require "correctness" with arbitrary parameters in $V_\kappa$ instead of just correctness of the theory $(V, \in, \kappa)$.

Just using a constant for the least reflective cardinal would not even be sufficient to define the satisfaction relation for ordinary set theory.

For more discussion, see (Taranovsky 2005) and also "Reflective Sequences" (especially "Introduction") below.

# 3 Axioms for Reflective Cardinals

We propose the following statements as axioms for reflective cardinals (see below for explanation and justification):

**A0.** ZFC
**A1.** $\forall x \exists y \forall z (z \in y \Leftrightarrow z \in x \wedge R(z))$
**A2.** $\forall a (R(a) \Rightarrow \mathrm{Ord}(a))$ [ that is $a$ is an ordinal ]
**A3.** $\forall a, b, c (\mathrm{Ord}(a) \wedge R(b) \wedge R(c) \wedge 0 < a < b < c \Rightarrow$
$(R(a) \Leftrightarrow \forall \hat{\varphi} \forall s \in V_a ((V_b \models \varphi(a, s)) \Leftrightarrow (V_c \models \varphi(b, s)))))$, where $\hat{\varphi}$ ranges over (codes for) formulas in set theory (without $R$) with all parameters shown.
**A4.** There is a proper class of reflective cardinals.
**A5.** Axiom schema of replacement for $R$-formulas.
**A6.** Axiom schema for every formula $\psi$ in the extended language (that includes $R$):
$\exists \kappa \in R \forall s \in V_\kappa (\psi(s) \Leftrightarrow V_\kappa \models \psi(s)))$ ($\kappa$ not free in $\psi$)
**A7.** Reflective cardinals are regular.
**A8.** Reflective cardinals are fully indescribable:
$\forall \kappa \in R \, \forall b \in R \backslash \kappa + 1 \, \forall \hat{\varphi} \, \forall A \subset \kappa \, V_b \models (\varphi(\kappa, A) \Rightarrow \exists \lambda < \kappa \, \varphi(\lambda, A \cap \lambda))$, where $\hat{\varphi}$ ranges over (codes for) formulas in set theory (without $R$) with all parameters shown.

**Notes:**
\* A6 ⇒ (A4 and A5); A8 ⇒ A7. All theorems will explicitly state which axioms (beyond ZFC) are used, and A4, A5, and A7 are included so that we can state theorems with less assumptions. A7 and A8 will be justified in the next section, where it will be shown that they do not increase the consistency strength of the theory.
\* Fully indescribable cardinals are totally indescribable. Full indescribability of κ can also be stated as a schema (implied by A8 for κ∈R): ($\varphi$ has two free variables and does not use $R$) $\forall A \subset \kappa \, (\varphi(\kappa, A) \Rightarrow \exists \lambda < \kappa \, \varphi(\lambda, A \cap \lambda))$. A fully indescribable cardinal need not be a $\Sigma_2$ elementary substructure of V: The least fully indescribable cardinal is below the least subtle cardinal (for each $n$, ZFC + subtle cardinal proves that the least $\Sigma_n^V$-indescribable cardinal is below the least subtle cardinal). However, if a fully indescribable cardinal is a $\Sigma_2$ elementary substructure of V, then it is strongly unfoldable.

A1 is obvious, and A2 is clear from the definition.
To explain A3, we start with the basic relations of reflective cardinals. Clearly,
(1) $R(a)$ and $R(b) \Rightarrow (\varphi(a) \Leftrightarrow \varphi(b))$ for every $\varphi$ in the language of set theory with parameters in $V_{\min(a,b)}$.
For the other direction, we want to formalize something like this:
(2) ($a < b$ and $R(b)$ and for all appropriate $\hat{\varphi}$ ($\varphi(a) \Leftrightarrow \varphi(b))) \Rightarrow R(a)$.
We create A3 by combining (1) and (2) and using a reflective cardinal above $b$ to formalize the quantification over formulas in (2). In A3, $\hat{\varphi}$ is the Gödel number of $\varphi$.
An alternative formulation of A3 is to use elementary substructures where $T(a, b)$ codes the theory of $V_b$ with constant $a$ and parameters in $V_a$:
$\forall a, b, c (\mathrm{Ord}(a) \wedge R(b) \wedge R(c) \wedge 0 < a < b < c \Rightarrow$
$(R(a) \Leftrightarrow (V_a, \in, T(a, b)) \prec (V_b, \in, T(b, c))))$.

A4 is the weakest existence axiom required for reflective cardinals to fulfill their role of giving semantics for higher order set theory (with arbitrary set-sized parameters). A4 (or a weakening stating that there is no largest reflective cardinal) is also necessary for A3 to imply (1) and (2).

ZFC together with A1-A4 forms a minimal theory of reflective cardinals, and this theory turns out to be finitely axiomatizable. However, the theory is missing replacement for $R$-formulas, and it does not even prove that some ordinal is a limit of reflective cardinals. In ZFC minus replacement, replacement is equivalent to this reflection schema: There is $V_a$ that is $\Sigma_n^V$ correct with parameters in $V_a$. A6 is the analogous schema for $R$-formulas, except that we also require $\kappa$ to be reflective. A6 is true since it is satisfied by every cardinal $\kappa$ with sufficiently strong reflection properties. We also include A5 in case the reader is not convinced about A6.

Arguments for large cardinal axioms often involve imagining $V$ as a completed totality, and then trying to go beyond and contemplating the largeness of the class of all ordinals. Since reflective cardinals provide semantics for higher order set theory, these arguments can be recast with greater precision as arguments about reflective cardinals. For example, using the semantics of higher order set theory, A7 is simply the statement that the axiom schema of replacement holds for arbitrary proper classes (in place of formulas). There is a general agreement that Ord should be regular (in so far as "Ord is regular" is meaningful), and inaccessible cardinals are already accepted by mainstream mathematicians in the form of Grothendieck universes.

## 4 Theory of Reflective Cardinals

The notion of reflective cardinals has a close connection with large cardinal axioms.

**Theorem 1** (ZFC + A1-A4): Let $\kappa$ be a reflective cardinal. Then,
**(a)** $V_\kappa$ satisfies ZFC, and (as a theorem schema) $(V_\kappa, \in)$ is an elementary substructure of $(V, \in)$.
**(b)** In HOD, $\kappa$ is fully indescribable.

**Proof: (a)** To see that $V_\kappa$ satisfies ZFC (as a single statement), note that by A3, if $\alpha < \beta < \gamma$ are reflective, then $V_\alpha$ is an elementary substructure of $V_\beta$. To see that $V_\kappa$ is $\Sigma_1$ correct, if a $\Sigma_1$ statement $\exists x \varphi(x)$ (possibly with parameter in $V_\kappa$) is true, then there is some $x$ witnessing its truthfulness and (by A4) two reflective cardinals above $x$, and so by A3, $\exists x \varphi(x)$ is true in $V_\kappa$. By induction on $n$, one sees (as a theorem schema) that $\kappa$ is $\Sigma_n$ correct for all $n$.
**(b)** Assume contrary, and let φ be a formula contradicting indescribability. Fix a definable well-ordering of HOD. For an ordinal $\alpha$ that is not φ-indescribable in HOD, let $A(\alpha)$ be the least witness to the describability of $A$; and let $A(\alpha)$ be 0 otherwise. If $a < b$ are reflective, then by definability of $A$, $A(a) = A(b) \cap a$. Thus, reflective cardinals are φ-indescribable with respect to $A$, and hence fully indescribable in HOD.

As the proof shows, reflective cardinals must satisfy large cardinal properties (such as indescribability) with respect to definable sets, and in HOD, that implies the corresponding large cardinal axioms. In general, however, the lifting from definable to arbitrary sets does not work, and ZFC + A0-A4 + no inaccessible cardinals appears to be consistent. I conjecture that the consistency strength of singular reflective cardinals is the same as the consistency strength of measurable reflective cardinals.
**Question:** What is the consistency strength of ZFC + A1-A5 + reflective cardinals are singular?

If we do not have a proper class of reflective cardinals, than besides A1-A3,A5, the axiomatization would also need to include a schema like A3 but using Ord and V in place of $c$ and $V_c$, and also a schema that for reflective $\kappa$, $V_\kappa$ is an elementary substructure of $V$. With this, a single reflective cardinal is (by compactness of first order logic) is equiconsistent with ZFC. Two reflective cardinals imply that reflective cardinals are (as a schema) fully indescribable in HOD.

**Theorem 2:** ZFC + A1-A4 is finitely axiomatizable.
**Proof:** For the axiomatization, we use A1-A4 plus (a finite axiomatization of) ZFC with separation and replacement limited to $\Sigma_1$ formulas. The inductive proof of Theorem 1(a) goes through in this theory, and thus every axiom of ZFC is provable.

**Theorem 3:** If $\kappa$ is a subtle cardinal, then for some regular $\lambda$ with $V_\lambda$ an elementary substructure of $V_\kappa$, there is $R \subset \lambda$ such that $(V_\lambda, \in, R)$ satisfies A1-A8.
**Proof:** Fix a well-ordering $W$ of $V_\kappa$ with sets of lower rank coming before sets of higher rank, and add a symbol for $W$ to the language. Let C = $\{\alpha < \kappa\colon V_\alpha$ is an elementary substructure of $V_\kappa$, allowing $W$ in formulas$\}$. Let $f(\alpha) \subset \alpha$ code the theory of $(V_\kappa, \in, W, \alpha)$ with parameters in $\alpha$, where $\alpha$ is an uncountable cardinal below $\kappa$. $C$ is closed and unbounded, so by Lemma 1.6 of "Subtle Cardinals and Linear Orderings", there is regular $\lambda$ in $C$ with a stationary subset $E$ of $\lambda \cap C$ homogeneous for $f$. Let $R$ be the following superset of $E$: $R = \{\alpha < \lambda\colon$ there is $\beta > \alpha$ in $E$ such that theory of $(V_\kappa, \in, \alpha)$ with parameters in $V_\alpha$ agrees with the corresponding theory for $\beta\}$. $(V_\lambda, \in, R)$ satisfies all the axioms. A6 holds because R is stationary, and A8 holds because every element of $V_\lambda$ is ordinal definable from $W$.

To get an exact equiconsistency, we introduce the following weakening of A3:
**A3a.** (Schema) $\forall a \in R \forall b \in R \forall s \in V_{\min(a,b)}(\varphi(a,s) \Leftrightarrow \varphi(b,s))$ ($\varphi$ has two free variables and does not use $R$).

"Ord is subtle" is the schema $\exists \alpha, b\, (a \in C, \beta \in C, a < b)\, f(V_\alpha) = f(V_\beta) \cap V_\alpha$ where $C$ and $f$ are definable (from parameters) classes and C is a club sublass of Ord and $\forall X\, f(X) \subset X$. It may also refer to the weaker schema $f(\alpha) = f(\beta) \cap \alpha$, which implies that for every set of ordinals $X$, $\mathrm{HOD}[X]$ satisfies the stronger schema (and hence both forms of ZFC + "Ord is subtle" prove the same $\Pi_2^V$ statements). We believe "Ord is subtle" should refer to the stronger schema; it is not even clear whether the weaker schema proves existence of $V_\kappa$ that satisfies ZFC.

**Theorem 4:** The following are equiconsistent and prove the same $\Pi_2^V$ statements:
(1) ZFC + A1 + A2 + A3a + 2 reflective cardinals + (schema) $R(\kappa) \Rightarrow (V_\kappa, \in) \prec (V, \in)$

(2) ZFC + A1 + A2 + A3a + A4-A7 + (schema) reflective cardinals are fully indescribable
(3) ZFC + (schema) Ord is subtle (either of the two forms)
(4) ZFC + (schema in $(V, \in, \kappa)$) $\kappa$ is fully indescribable + (schema) $V_\kappa \prec V$

**Proof:** Let (3) use the stronger form of Ord is subtle. We have (2)⇒((1) and (3) and (4)). (4)⇒(3). For each $n$, (3) proves existence of two cardinals satisfying (1) restricted to $\Sigma_n$ formulas, and thus (3)⇒(1) for statements in $(V, \in)$ (that is statements in $(V, \in)$ provable in (1) are provable in (3)). Let $X$ code a well-ordering of $V_\alpha$ where $\alpha$ is the least ordinal with $V_\alpha \prec_{\Sigma_2} V$. (1) implies that $\text{HOD}[X]$ satisfies (4) (with reflective cardinals being fully indescribable), and hence (1)⇒(4) for $\Pi_2^V$ statements (see also proof of Theorem 5). Finally, given (4), for a natural number $n$, if we define $R$ below $\kappa$ based on A3a restricted to $\Sigma_n$ formulas, then $R$ is stationary below $\kappa$ and $V_\kappa$ satisfies (2) restricted to $\Sigma_n^R$ statements. Thus (4) ⇒ (2) for statements in $(V, \in)$, which completes the proof.

**Notes:**
* (2) and (4) prove the same $(V, \in, \kappa)$ statements where for (2) $\kappa$ denotes a reflective cardinal. Similarly, (1) and ZFC + A1 + A2 + A3a + A4-A6 prove the same $(V, \in, \kappa)$ statements.
* (1) and (3) (stronger form) prove the same (V, ∈) statements. (Given $f$ and $C$, any pair of reflective cardinals above the parameters used to define $f$ and $C$ witnesses "Ord is subtle" for $f$ and $C$, so (1)⇒(3).)
* (2) can be strengthened with the schema: there is $\Sigma_n^{V,R}$-correct reflective cardinal $\kappa$ such that ($\varphi$ has two free variables and does not use $R$)
$\forall A \subset \kappa \, (\varphi(\kappa, A) \Rightarrow \exists \lambda < \kappa \, (R(\lambda) \land \varphi(\lambda, A \cap \lambda)))$ without introducing new $(V, \in, \kappa)$ theorems.
* If the requirement of being an elementary substructure of V is dropped from (1) and (4), then the two theories remain equiconsistent with each other, and we conjecture remain equiconsistent with the other theories.

**Theorem 5:** In each of ZFC+A1-A4 and ZFC+A1-A5 and ZFC+A1-A6, axiom A8 (reflective cardinals are fully indescribable) is $\Pi_2^V$ (for formulas without $R$) conservative over the base theory.

**Proof:** Let $T$ be a $\Pi_2^V$ statement that is unprovable in the base theory, and let $M$ be a model in which $T$ fails. Let $X$ — as computed in $M$ — be a set of ordinals coding a well-ordering of some $V_\kappa$ that witnesses the failure of $T$. Then $(\text{HOD}[X]^M, \in, R')$ satisfies A8 and the failure of T. Here $R'$ is defined in $M$ as the natural conversion of $R$ to HOD[X]: $R'(a) \Leftrightarrow$ there is $b,c$ such that $R(b) \land R(c) \land c > b > \max(a, \sup(X))$ and as computed in $\text{HOD}[X]^M$, the theory of $(V_c, \in, a)$ with parameters in $V_a$ agrees with that of $(V_c, \in, b)$.

$\Pi_2^V$ is a very broad class of statements (it even includes GCH), so Theorem 5 suggests that A8 is "safe". The known arguments for large cardinal axioms (see for example [Kanamori 2008]) combined with Theorems 1-5 create a powerful argument for acceptance of totally indescribable (and even strongly unfoldable) cardinals as axioms. Being totally indescribable corresponds to our intuition about the theory of finite types above $V$, and without V=HOD, axioms A1-A6 do not be appear to imply any boldface indescribability.

Once a large cardinal notion is accepted, the use of reflective cardinals can often be used to argue for a stronger notion, as illustrated by the next theorem.

**Theorem 6** (ZFC+A1-A4): **(a)** If $j: V \to M$ is elementary with critical point $\kappa$ and $R(\kappa)$ and $j(R)(\kappa)$, then $o(\kappa) > (\kappa_{\text{plus}}^+)^L$ where $\kappa_{\text{plus}} = (\kappa^+)^{\text{HOD}}$.
**(b)** If V=HOD, and a reflective cardinal is $R$-strong, then reflective cardinals are Shelah (and hence Woodin).
**Proof: (a)** Assume contrary. Because $0^\#$ exists, ordinals below $(\kappa_{\text{plus}}^+)^L$ can be coded (in a definable way absolute between V and M) by ordinals below $\kappa_{\text{plus}}$. In turn, an ordinal below ordinal below $\kappa_{\text{plus}}$ can be represented by an ordinal definable well-ordering of $\kappa$ of the right order type (note that since $\kappa$ is the critical point, $\text{HOD}^M \cap P(\kappa) = \text{HOD} \cap P(\kappa)$). Let $A(\kappa)$ be such a representation of $o(\kappa)$. We have $A(\kappa) = j(A(\kappa)) \cap \kappa = A^M(j(\kappa)) \cap \kappa = A^M(\kappa)$ (the last equality holds because both $\kappa$ and $j(\kappa)$ are reflective in $M$), which contradicts $o(\kappa)^M < o(\kappa)$.
**(b)** Let $\kappa, \lambda, \mu$ be reflective, $\kappa < \lambda < \mu$, and $j: V \to M$ a $\mu$-$R$-strong elementary embedding with crit($j$) = $\kappa$. Let $f_a : a \to a$ be the least (under the canonical well-ordering of V (V=HOD)) witness (if any) that $a$ is not Shelah. Since $j$ is $\mu$-$R$-strong, $f_\kappa^M = f_\kappa$ and $f_\lambda^M = f_\lambda$. Furthermore, $j(f_\kappa) \cap V_\lambda = f_\lambda$, so $j(f_\kappa)(\kappa) = f_\lambda(\kappa) < \mu$ as required.

**Note:** Because ZFC+A1-A4 lacks full separation, it does not prove that a $\kappa$-complete ultrafilter generates an embedding that is elementary for arbitrary $R$-formulas (where $j(R)$ is used in $M$ in place of $R$). However, the theorem stays true if we only require $j$ to be elementary for formulas without $R$.

## 5 Iterations of Reflectiveness

### 5.1 Axiomatization and Properties

The same process that led us to reflective cardinals can be repeated any number of times. For example, we can call a cardinal $\kappa$ is 2-reflective iff $(V, \kappa, \in, R)$ has the same theory with parameters in $V_\kappa$ as $(V, \lambda, \in, R)$ for every cardinal $\lambda > \kappa$ with sufficiently strong reflection properties ($R$ is the predicate for reflective cardinals). Here is a basic axiomatization of arbitrary ordinal levels of reflectiveness:

**B1.** ZFC
**B2.** for all $x$ there is $\{(y, z) \in x : R(y, z)\}$
**B3.** $R(y, z) \Rightarrow \text{Ord}(y) \wedge \text{Ord}(z)$
**B4.** $R(0, x)$ iff $x$ is a regular uncountable cardinal
**B5.** $R(n+1, a) \Rightarrow R(n, a)$
**B6.** for limit $n$, $R(n, a) \Leftrightarrow \forall m < n R(m, a)$.

**B7.** $\forall n, a, b, c (\mathrm{Ord}(n) \wedge \mathrm{Ord}(a) \wedge R(n+1, b) \wedge R(n+1, c) \wedge 0 < a < b < c \Rightarrow (R(n+1, a) \Leftrightarrow \forall \hat{\varphi} \forall s \in V_a((( V_b, \in, R_n) \models \varphi(a, s)) \Leftrightarrow ((V_c, \in, R_n) \models \varphi(b, s)))))$
where $\hat{\varphi}$ ranges over (codes for) formulas in set theory with $R_n$ with all parameters shown.

**B8.** Schema ($m$ is a natural number): There is a $\Sigma_m^{V,R}$ correct $\kappa$-reflective cardinal $\kappa$.

**Notes:**
* $R(y, z)$ means that $z$ is level $y$-reflective. If $R_y(z)$ is $R(y, z)$, then $R_1$ is the predicate for reflectiveness used in the previous section.
* B7 is the analogue of A3.
* B8 is true because it is satisfied by every ordinal with sufficiently strong reflection properties.
* For a predicate $P$, $\kappa$ is defined to be $\Sigma_n^{V,P}$ correct iff for all $\Sigma_n^{V,P}$ $\varphi$ (with one free variable) for all $s \in V_\kappa$, $(V_\kappa, \in, P)$ satisfies $\varphi(s)$ iff $\varphi(s)$ holds in $V$. $\Sigma_n^{V,P}$ formulas are formulas of the form $\exists x_1 \forall x_2 \ldots x_n \psi$ ($n$ alternating quantifiers) where $\psi$ is a bounded quantifier formula (but allowing $P$).
* A weakening of B8 that would allow finite axiomatizability is:
$\forall x \in \mathrm{Ord} \, \exists y > x \, R(y, y)$.

Axioms B1-B8 seem too incomplete. For example, they do not appear to answer whether reflective cardinals are Mahlo, or whether there are stationary many reflective cardinals. To resolve the incompleteness, we add the following strengthening of B8:

**B9.** Schema ($n$ is a natural number): There is a $\Sigma_n^{V,R}$-correct $\kappa$-reflective measurable cardinal $\kappa$, and a corresponding embedding $j$ (with crit($j$)=$\kappa$) such that $j(R)(\kappa, \kappa)$ holds.

Because of its correctness, $\kappa$ in B9 can represent the class of ordinals Ord, and $j(R)(\kappa, \kappa)$ captures the idea that reflective cardinals are those that are most like Ord.

B1-B8 plus any one of B9 or B9a (see below) or V=HOD[X] implies that $n+1$-reflective cardinals are fully indescribable, even for formulas in $(V, \in, R_n)$.

We start our investigation with finite levels of reflectiveness. For finite $n$, there is a natural weakening of reflectiveness: $\kappa$ is weakly $n$-reflective cardinal iff the theory of $(V, \in, \kappa, a_1 \ldots a_{n-1})$ with parameters in $V_\kappa$ is correct (as defined using $n$-reflective cardinals) where $a_1 < a_2 < \ldots a_{n-1}$ are $n-1$ reflective cardinals above $\kappa$. An axiomatization of weakly $n$-reflective cardinals uses A0-A2, A4, and a natural modification of B7:

**B7a:** $R(a)$ holds iff for all $a_1 \ldots a_{n+1}$ with $R(a_i)$ and $a \leq a_1 < a_2 < \ldots a_{n+1}$, the theory of $(V_{a_{n+1}}, \in, a, a_2 .. a_n)$ with parameters in $V_a$ agrees with the theory of $(V_{a_{n+1}}, \in, a_1, a_2 .. a_n)$.

$n+1$-weakly reflective cardinals are n-reflective, and $n$-reflective cardinals are weakly $n$-reflective. While for $n > 1$, the notion of $n$-reflective is more expressive than that of weakly $n$-reflective, the difference is slight, much less than the difference between $n$-reflective and $n + 1$-reflective. To see that, if we extend expressiveness slightly by adding satisfaction relation for $(V, \in)$, then the extended weakly $n$-reflective cardinals become more expressive than the (unextended) $n$-reflective cardinals. By induction on $m$, one sees that an $m$-tuple of extended weakly n-reflecting cardinals ($m \leq n$) can be used to define the predicate for $m$-reflective cardinals below it, and with the least element of the $m$-tuple corresponding to an elementary substructure of $(V, \in, R_m)$. I am not sure whether weakly $n$-reflective cardinals should simply be called $n$-reflective, but this paper currently uses $n$-reflective for the stronger notion.

The hierarchy of $n$-reflective cardinals is intertwined with the hierarchy of subtle/ineffable cardinals.

**Theorem 7** (A0-A2, A4, B7a(n+1)): weakly $n + 1$-reflective cardinals are $n$-ineffable in HOD.
**Proof:** Fix a definable well-ordering $W$ of HOD. We first show that the reflective cardinals are $n$-subtle. Let $C(\alpha)$ and $f(\alpha)$ witness least failure of $n$-subtlety of $\alpha$. Fix $a_1 < a_2 < \ldots < a_{n+2}$ that are $n + 1$-weakly-reflective. By reflectiveness (and definability and closure and unboundness of $C$), $a_1 \ldots a_{n+1}$ belong to $C(a_{n+1})$. By weak $n + 1$-reflectiveness of $a_i$, $f(a_{n+2})(a_1 \ldots a_n) = a_1 \cap f(a_{n+2})(a_2 \ldots a_{n+1})$, which shows that $a_{n+2}$ is $n$-subtle.
To get ineffability, we need stationary many cardinals that are sufficiently reflective cardinals, so we introduce the following weakening: $R'(a)$ iff the theory of $(V_{a_n+\omega}, \in, a, a_2 \ldots a_n)$ with parameters in $V_a$ is correct (that is agrees with the theory for weakly $n$-reflective $a$). The proof is completed by checking that $R'$ is stationary below weakly $n$-reflective cardinals and that the $R'$ cardinals suffice for the agreement of subsets for ineffability.

**Theorem 8:** If $\kappa$ is an n-subtle cardinal, then for some regular $\lambda$ with $V_\lambda$ an elementary substructure of $V_\kappa$, there is $R \subset \lambda$ such that $(V_\lambda, \in, R)$ satisfies A0-A2, A6, B7a. Alternatively, $R$ can be chosen to satisfy B1-B8, restricted to $n$-reflective cardinals or less. Moreover, the $n$-reflective (and weakly $n$-reflective) cardinals can be chosen to be $n - 1$-ineffable (if $n > 1$) and fully indescribable in $V_\lambda$.
**Proof:** Analogous to the proof for consistency of reflective cardinals from a subtle cardinal. For $n$-reflective cardinals, use their reduction to a slightly stronger variant of weakly $n$-reflective cardinals. To get ineffability, use a well-ordering $W$ and the proof of Theorem 7 (and analogously with the indescribability).

A substantial fragment of the theory of $n$-reflective cardinals is conservative over ZFC:
**Theorem 9:** The following theory is conservative over ZFC:
- The language consists of first order logic, $\in$, and $\omega$ predicates $R_1, R_2, \ldots$.
- ZFC, plus replacement for formulas in the extended language
- (schema, $n > 0$) there is cardinal κ satisfying $R_n(\kappa)$
- (schema, $n > 0$) $\forall \kappa \in R_{n+1}\ \kappa \in R_n$
- (schema, $n > 0$, and '$\prec$' uses a schema)
$\forall \kappa \in R_n\ (V_\kappa, \in, R_1, \ldots, R_{n-1}) \prec (V, \in, R_1, \ldots, R_{n-1})$.

- (schema, $n \in \mathbb{N}$, $\varphi$ is a formula in $(V, \in, R_1, \ldots, R_n)$ with two free variables, and $S$ ranges over sets definable in $(V, \in, R_1, \ldots, R_n)$)
$\forall \kappa \in R_{n+1} \forall \lambda \in R_{n+1} \forall s \in S\, (\varphi(\kappa, s) \Leftrightarrow \varphi(\lambda, s))$.

**Proof:** Given a model of ZFC and a finite fragment of the theory, we interpret the fragment as follows, and the conservation result follows. Let $n \in \mathbb{N}$ be sufficiently large relative to the fragment, $m \in \mathbb{N}$ the largest integer such that $R_m$ is used in the fragment, and $\kappa$ an ordinal such that $(V_\kappa, \in) \prec_{\Sigma_{n(m+1)}} (V, \in)$. Let T be the $\Sigma_n$ theory $(V, \in, \kappa)$ with parameters in $V_\gamma$ where $\gamma$ is the least ordinal such that $(V_\gamma, \in) \prec_{\Sigma_n} (V, \in)$ and set $R_1(\alpha)$ to hold iff $\Sigma_n$ theory of $(V, \in, \alpha)$ with parameters in $V_\gamma$ is T. Define $R_{i+1}$ analogously, always using $V, \in, R_1, \ldots, R_i$ in place of $V, \in$, and verify that the fragment of the theory holds.

**Lightface Reflective Cardinals**

We briefly discuss what happens to the consequences provable from the axioms if in the definition of reflective cardinals, we simply used 'theory' instead of 'theory with parameters'.

* As the theorem above shows, a version of the theory (using $\mathbb{N}$ in place of $S$) of lightface $n$-reflective cardinals is conservative over ZFC. However, to make the theory appear more definite, we would want to replace that schema with the appropriate analogue of B7 that does not use $s$ as the second argument of $\varphi$. The consistency strength remains below inaccessible. If $\rho$ is inaccessible in $L_{\omega+2}(V_\rho)$, then we can pick $\kappa < \rho$ such that $V_\kappa \prec_{\Sigma_{\omega^2+\omega}} V_\rho$ and set $R_n(\alpha)$ based on agreement of $\Sigma_{\omega n}$ theory of $(V_\rho, \in, \alpha)$ with $(V_\rho, \in, \kappa)$. The extension to transfinite levels of lightface reflectiveness is straightforward. If we require lightface reflective cardinals to be regular, then it suffices for $\rho$ to be Mahlo in $L_{\omega+2}(V_\rho)$ (so the consistency strength remains below Mahlo).

* For $n > 1$, there is a strengthening of lightface $n$-reflectiveness that can be axiomatized as follows: (schema; $\varphi$ has two free variables and does not use $R$) $\forall \kappa \in R_n^n \forall \lambda \in R_n^n \forall m \in \mathbb{N}\, (\varphi(\kappa, m) \Leftrightarrow \varphi(\lambda, m))$, where $R_n^n$ is the predicate for $n$-element subsets of $R_n$, or a single statement that encompasses the schema (assuming $\forall \kappa \in R_n\, (V_\kappa, \in) \prec (V, \in)$) by using $V_\mu$ in place of $V$ for a lightface reflective $\mu > \max(\kappa \cup \lambda)$. A model of this ($n < \omega$ arbitrary but fixed; the model satisfies ZFC + a proper class of such cardinals + $\forall \kappa \in R_n\, (V_\kappa, \in) \prec (V, \in)$, etc.) can be obtained using a weakly compact cardinal. Conversely, under V=HOD, given the schema (even without $\forall \kappa \in R_n\, (V_\kappa, \in) \prec (V, \in)$) and a lightface 3-reflective limit $\kappa$ of $\kappa$ lightface 3-reflective cardinals, lightface reflective cardinals are weakly compact (because they would form a counterexample to the least failure of the partition relation). Furthermore, under V=HOD and the schema, if the set of lightface $n+2$-reflective cardinals is stationary below (a regular lightface $n+2$-reflective) $\kappa$, then they are $n$-ineffable. The arguments below about $\alpha$-Erdös cardinals (and $0^\#$) also apply to lightface $\omega$-reflective cardinals, assuming that they are required to be indiscernibles for predicates definable in $(V, \in)$ without parameters.

## 5.2 $\omega$-Reflective Cardinals and Beyond

Existence of a proper class (or just $\omega_1$) of $\omega$-reflective cardinals is inconsistent with $V = L$.

**Theorem 10** (B1-B7 + a proper class of $\omega$-reflective cardinals + the full axiom of replacement): Every set has a sharp. For every countable-in-HOD ordinal $\alpha$, there is $\alpha$-Erdös cardinal in HOD. If V=HOD and there is a reflective limit of $\omega$-reflective cardinals, then reflective cardinals are Ramsey.
**Proof:** For a set of ordinals $X$, $\omega$-reflective cardinals above $X$ (that is above $\sup(X) + 1$) form indiscernibles for $L[X]$, so $X$ has a sharp.
If there is no $\alpha$ Erdős cardinal in HOD, let $A(\kappa)$ be the least witness of the failure of the partition relation. Then one easily shows that for a $\omega$-reflective $\kappa$, $\omega$-reflective cardinals below $\alpha$ are indiscernibles for $A(\kappa)$, and thus the partition relation holds in $V$. Since $\alpha$ is countable in HOD, it also holds in HOD and hence there is an $\alpha$-Erdős cardinal in HOD.
If $\kappa$ is not Ramsey in HOD, let $A(\kappa)$ be the well-ordering-least partition relation contradicting Ramseyness. The $\omega$-reflective cardinals below $\kappa$ then form the desired homogeneous set.

For the proof of existence $\alpha$-Erdös cardinal in HOD (for $\alpha$ countable-in-HOD), it is enough to have $\alpha + \omega$ $\omega$-reflective cardinals.
We do not know if "B1-B7 + a proper class of $\omega$-reflective cardinals + reflective limit of $\omega$-reflective cardinals" proves that reflective cardinals are Ramsey in HOD.

**Theorem 11:** If κ is $\geq \alpha$-Erdös cardinal ($\alpha$ a limit ordinal), then there is $R$ such that $V_\kappa$ satisfies B1-B7 + there are $\alpha$ $\omega$-reflective cardinals.
**Proof:** Choose a well-ordering $W$ of $V_\kappa$ with sets of lower rank coming before sets of higher rank. Let $f$ code theory of $(V_\kappa, \in, W)$ with finite sets of ordinals as parameters. Get a homogeneous set $E$ for $f$ of order type $\alpha$, and let $M$ be the Skolem Hull of that set. Every order preserving one-to-one mapping $E \to E$ leads to an embedding of $M$ into $M$, so let $E'$ be the set of critical points of such embeddings. Then all members of $E'$ satisfy the requirements for ω-reflective cardinals.

**Corollary 12:** If $\kappa$ is Ramsey, then there is $R$ such that $V_\kappa$ satisfies ZFC+B1-B7 + there is a proper class of $\omega$-reflective cardinals.

Natural semantics for the axioms is provided by measurable cardinals:

**Theorem 13:** If $\kappa$ is measurable, then there is $R$ such that $(V_\kappa, \in, R)$ satisfies B1-B8.
Moreover, if $j : V \to M$ is an elementary embedding with critical point $\kappa$, then there is a unique (as a predicate on $V_\kappa$) $R$ satisfying B1-B8 and $j(R)(\kappa, \kappa)$.
**Proof:** By induction on degree of reflectiveness. Suppose that $R_n$ is uniquely defined and $j(R_n)(\kappa)$ holds. Then in $M$, let $\kappa$ be $n+1$-reflective. Using $V^M_{j(\kappa)}$ and $j(R_n)$, this uniquely defines $R_{n+1}$ below $\kappa$; moreover $R_{n+1}$ has measure 1, so $j(R_{n+1})(\kappa)$ holds.

Following Mitchell, a normal measure $U$ on $\kappa$ is called a weak repeat point iff for every $X \in U$ there is a normal measure $W$ on $\kappa$ in $\text{Ult}(V, U)$ such that $X \in W$.

**Theorem 14:** If $\kappa$ has a normal measure that is a weak repeat point, then there is $R$ such that $(V_\kappa, \in, R)$ satisfies B1-B9.
**Proof:** Apply Theorem 13 to define $R$ from the measure. The set
$\{x < \kappa : R(x, x) \land (V_x, \in, R|x) \prec (V_\kappa, \in, R)\}$ has measure one, so in $M$ there is an ultrafilter on $\kappa$ that satisfies B9.

The notion of α-reflective cardinals is more forward looking than simply the notion of reflective cardinals, and given that for α≥ω it is necessarily stronger than $0^\#$, the use of measurable cardinals in B9 is not a dramatic strengthening. However, there is a weakening of B9 that does not require existence of measurable cardinals.
**B9a:** Axiom schema consisting of (or equivalent to) all theorems provable in ZFC about $(V_\kappa, \in, R)$ where $R$ satisfies B1-B8, and there is an elementary embedding $j : V \to M$ with crit(j)=κ and $j(R)(\kappa)$.
**Note:** We can further weaken B9a by using ZFC\P (or other appropriate theory) in place of ZFC.
B9a is equiconsistent with ZFC + measurable cardinal, but it has a somewhat unsatisfying feel of using a measurable cardinal without claiming existence of a measurable cardinal. That is why we chose the stronger B9.

## 5.3 HOD Sharp

Although I defined reflective cardinals in $V$, they (or their analogues) can also be defined for canonical inner models, including the constructive universe $L$. If $M$ is is a canonical inner model (such as $L$), then the indiscernibles derived from $M^\#$ will be ω-reflective in $M$. For each finite level $n$, the predicate for $n$-reflectiveness can be added to $M$ in the unique way that satisfies the axioms and agrees with the indiscernibles; and M satisfies replacement for formulas involving $R_n^M$). However, for many $M$ (including $M = L$, but not $K^{\text{DJ}}$) the predicate for ω-reflective ordinals would allow us to transcend $M$. (For $M = L[U]$, ω-reflective ordinals transcend $M$, but $K^{\text{DJ}}$ is obtained from $L[U]$ by iterating away $U$, and the missing large cardinal structure is the reason ω-reflective ordinals for $K^{\text{DJ}}$ do not lead to new sets.)

It is plausible that ordinal definable real numbers are precisely those that are first order definable using a finite level of reflectiveness. If this is so, then we define $\text{HOD}^\#$ (HOD sharp) to be the real number coding
$\{\hat{\varphi} : n \in \mathbb{N} \land (V_k \models \varphi(k_1, \ldots, k_n))\}$ where $k_1 < k_2 < \ldots < k_n < k$ are $n$-reflective $\}$.
$\text{HOD}^\#$ is the canonical real number that codes all ordinal definable real numbers. $\text{HOD}^\#$ is not ordinal definable, and our apparent ability to define $\text{HOD}^\#$ illustrates the power of our extension to the language of set theory.

To extend the notion of $\text{HOD}^\#$ to arbitrary sets, we propose the following conjecture, the consistency of which is an open problem.
**Conjecture:** For every nonzero ordinal $\kappa$, every ordinal definable subset of $\kappa$ is definable (in $(V, \in)$) from $\kappa$, an element of $\kappa$, and a finite set of $\omega$-reflective cardinals above $\kappa$.

The conjecture implies that for every infinite $\kappa$, $\text{cf}((\kappa^+)^{\text{HOD}}) = \omega$. Assuming GCH in HOD (which is likely to be true), the conjecture is equivalent to $(\kappa^+)^{\text{HOD}} = \sup\{$ $\lambda$: a bijection between $\kappa$ and $\lambda$ is first-order-definable in $(V, \in)$ using a finite set of $\omega$-reflective cardinals above $\kappa\}$.

A strengthened version of the conjecture is that for every nonempty set of ordinals $S$, every ordinal definable from $S$ subset of $\sup(S)$ is definable (in $(V, \in)$) from $S$, an element of $\sup(S)$, and a finite set of $\omega$-reflective cardinals above $\sup(S)$.

Another way to propose $\text{HOD}^{\#}$ is to claim a nontrivial elementary embedding HOD→HOD, which is expressible in the higher order set theory. However, it is not clear whether we can actually define such an embedding in the extended language (and it is known not to be ordinal definable).

The notion of reflective cardinals (and the iterations of this notion) is the most expressive extension to the language of set theory for which a reasonable axiomatization is known. There is of course incompleteness in the axiomatization, but it is comparable in scope to what we have in the rest of set theory.

# 6 Reflective Cardinals and Supercompactness

Reflective cardinals can be axiomatized at any large cardinal level beyond the intrinsic strength of the basic axioms. By reflection, they satisfy all large cardinal properties that are fully realized in $V$. In this section, we use reflective cardinals to argue for existence supercompact cardinals, and propose an axiomatization of reflective cardinals using supercompact cardinals. I should note that the argument is not conclusive. Supercompact cardinals probably exist, but we are not sure.
**Note:** Reinhardt's argument for extendible cardinals and Vopenka's argument for Vopenka's cardinals can be viewed as precursors to this argument.

Let $\lambda$ be a reflective cardinal, and $\beta$ be a reflective (or just $\Sigma_2^V$ sound) cardinal above $\lambda$. By largeness of $\lambda$ and by its reflecting properties, there must be a cardinal $\kappa < \lambda$ that is very similar to $\lambda$—as similar as it is possible to be and as is expressible in the language of set theory. (Technically, "expressible in the language of set theory" should be replaced here with something like "expressible with $\Sigma_{100}^V$ formulas".) We want to formalize that similarity as much possible.
By similarity of $\kappa$ and $\lambda$, and by the symmetry of $V$, for some $\alpha > \kappa$ (and with $\alpha < \lambda$), the structure of $(V_\alpha, \in, \kappa)$ is very similar to that of $(V_\beta, \in, \lambda)$. Similarity between structures is best expressed by using mappings between structures. Here, that means an elementary embedding $j : (V_\alpha, \in, \kappa) \to (V_\beta, \in, \lambda)$. By maximality of $V$, $j$ should exist for some $\alpha$ if $(V_\alpha, \in, \kappa)$ can be sufficiently similar $(V_\beta, \in, \lambda)$.
The embedding $j$ must have a critical point, $\text{crit}(j)$. The natural choice of $\text{crit}(j)$ is $\kappa$, and by the reflecting properties of $\kappa$ this choice should be possible (for appropriate $\alpha$).

Moreover, having $j$ with crit($j$)=$\kappa$ corresponds to the semantics of higher order set theory, where sets stay the same while realizations of proper classes and higher level structures depend on the choice of the reflective cardinal. If Ord (the class of all ordinals) is represented by $\kappa$, then $j$ is just a change of representation by which Ord becomes represented by $\lambda$, and by which proper classes and beyond are moved from being sets in $V_\alpha$ to being sets in $V_\beta$.

Thus (in the absence of some unknown problem), we have non-trivial elementary $j : V_\alpha \to V_\beta$ with $\lambda = j(\mathrm{crit}(j))$ and with $\beta$ being $\Sigma_2^V$ sound. This means that $\lambda$ (and thus every reflecting cardinal) is supercompact.

**Note:** Because $\kappa \neq \lambda$, there must be a limit to the degree of similarity between them. We stay within the limit because we are only requiring only that (essentially) there is a structure of the form $V_a$ that is similar enough to $V_\beta$, and not requiring that the similarity (and $j$ in particular) can be extended arbitrarily high up the cumulative hierarchy.

Existence of a reflective supercompact cardinal is a much broader assertion than that of existence of a supercompact cardinal. For example, existence of a supercompact cardinal does not imply existence of a model of ZFC that satisfies all true $\Sigma_3^V$ statements. (Analogously, Con(ZFC + supercompact cardinal) does not imply existence of a transitive model of ZFC.) Even having a reflective supercompact cardinal is not the highest level of expressiveness to which supercompactness can be added.

For a predicate $P$, let us say that a cardinal $\kappa$ is $P$-supercompact iff for every $l > \kappa$, there is $a$ and elementary j: $(V_a, \in, P) \to (V_l, \in, P)$ with $\kappa = j(\mathrm{crit}(j))$. The argument for existence of supercompact cardinals can be modified to show existence of $R$-supercompact cardinals.

To extend the strength of supercompactness to the full expressiveness of the language, we strengthen A5 to read:

**P1:** Schema ($n$ is a natural number): There is $\Sigma_n^{V,R}$ correct $R$-supercompact $\kappa$ such that $R(\kappa)$ holds.

P1 is probably true, but we do not have enough confidence to add it as axiom: We are not yet able to rule out some kind of impossibility that would prohibit the required elementary embeddings.

P1 actually goes slightly beyond supercompactness; for example, all $R$-supercompact cardinals are extendible.
An important consequence of P1 is that for a proper class of reflective cardinals, for every $a < b < c$ with $a$ belonging to that class and $c$ reflective, there is $d > c$ and elementary $j : V_b \to V_d$ with $\mathrm{crit}(j) = a$ and $j(a) = c$.

Existence of an $R$-supercompact reflective cardinal implies Vopenka's principle for definable proper classes — more precisely, definable at a lower complexity level than what is needed to define $R$. However, existence of a Vopenka cardinal is a philosophically stronger assumption — it assumes that elementary embeddings can be found from just the largeness of the cardinal without relying on the symmetry of the structure. The axiom of choice is notorious for breaking symmetry, so existence of a Vopenka cardinal does not immediately follow from Vopenka's principle for definable classes. However, I am not suggesting that Vopenka cardinals do not exist. For example, if V=HOD, then reflective cardinals are Vopenka (assuming existence of $R$-supercompact reflective cardinals).

As I wrote in section "Iterations of Reflectiveness", reflectiveness can be iterated any number of times. The analogue of B9 is:
**P2:** There is $\Sigma_n^{V,R}$ correct $R$-supercompact $\kappa$ such that $R(\kappa, \kappa)$ holds.

# 7 Beyond reflective cardinals

## 7.1 Overview

By the undefinability of truth, whenever we define a valid formal language, we can go beyond it. However, at the outer limits of current mathematical expressiveness, we have the notions but not a reasonable axiomatization. It is these notions that I present here, hoping that they will be axiomatized in the future. (Update: A partial axiomatization is in "Reflective Sequences" below.) The key idea is that instead of using a single reflective cardinal, we use a sequence of them. Since we do not have a good axiomatization, some of the discussion here is philosophical and may appear vague to the reader, but the theorems are of course unambiguous.

To motivate the extensions, here is a definition that is expressible using iterated reflectiveness:
**Definition**: A set $S$ is a *reflective n-tuple* iff there are $n$-reflective $b_1, \ldots, b_{n+1}$ with $b_1 < b_2 < \ldots b_{n+1}$ and such that $(V_{b_{n+1}}, \in, S)$ with parameters in $V_{\min}(S)$ agrees with the theory in $(V_{b_{n+1}}, \in, T)$ where $T$ is $\{b_1, \ldots, b_n\}$, and for every ordinal $\lambda$, $S \backslash \lambda$ is a reflective tuple under the above definition.

**Notes:**
\* Reflective $n$-tuples correspond to weakly $n$-reflective cardinals. Thus, as discussed in a previous section, the predicate for reflective $n$-tuples is slightly less expressive (when combined with ordinary set theory) than the predicate for $n$-reflective cardinals. However, even a single reflective n+1-tuple S can be used to define the true predicate for $n$-reflective cardinals for ordinals below min($S$).
\* The requirement for reflectiveness of $S \backslash \lambda$ does not change the expressive power of the predicate but allows us to have a nicer theory.

In general, given a type of objects, the reflection condition is formulated in the following way:
$R(\mathbb{O})$ iff $\mathbb{O}$ satisfies sufficiently strong reflection properties so that the theory $(V, \in, \mathbb{O})$ with parameters in $V_{\min}(\mathbb{O})$ is independent of $\mathbb{O}$ (for objects of the required type).
(Note: For technical convenience, we may also want to require that for every ordinal $\lambda$, the theory $(V, \in, \mathbb{O} \backslash V_\lambda)$ with parameters in $V_\lambda$ is correct for the type of $\mathbb{O} \backslash V_\lambda$.)

An axiomatization would consist of:

1. Description of the type of $\mathbb{O}$ for $\mathbb{O}$ satisfying $R$ (for example $\mathbb{O}$ is an ordinal).
2. The reflection condition above (analogous to A3).
3. The analogue of A5 for $R$ (such as replacement for $R$-formulas).
4. Large cardinal and reflection properties for $\mathbb{O}$ satisfying $R$.
5. Existence axioms for $\mathbb{O}$ satisfying $R$.

(1), (2), and (3) are routine. The problem is that we are not sure about (4) and (5).

## 7.2 Reflective $\omega$-Sequences

The first notion, and the one that we understand best, is
**Notion 1:** A reflective sequence of ordinals of length $\omega$.

This notion is more expressive than any ordinal iteration of reflectiveness: If $S$ is a reflective sequence of length $\omega$, then the two variable $R$ (in section "Iterations of Reflectiveness") can be defined from $S$ for ordinals below $\min(S)$ (and under our definition below $\sup(S)$).
To axiomatize such sequences, for (4), one has to specify not merely the properties with the individual members of $S$, but also the properties of $S$ as a whole. For example, if there are I3 embeddings, then likely $S$ is the critical sequence for some elementary $j: V_\lambda \to V_\lambda$. On a more prosaic level, $S$ should be Prikry generic for some $L[U]$ where $U$ is a measure in $L[U]$. In fact, it appears that a basic theory requires little more than a measurable cardinal.

**Definition:** For a set $A$, a set of ordinals $S$ of order type $\omega$ is *lightface reflective for* $A$ iff the theory $(V, \in, S, A)$ agrees with the theory $(V, \in, T, A)$ for every set of ordinals $T$ of order type $\omega$ with sufficiently strong reflection properties. $S$ (of order type $\omega$) is *reflective* iff for every $A \in V_\alpha$ with $\alpha < \sup(S)$, $S \backslash \alpha$ is lightface reflective for $A$.

The definition does not immediately imply much about the structure theory for such sequences. Thus, we start with a conjecture.
**Conjecture:** A set of ordinals of order type $\omega$ with each member belonging to a reflective $\omega$-sequence is a reflective $\omega$-sequence.

As the theorem below shows, any natural theory of reflective $\omega$-sequences is inconsistent with V=HOD. However, the concept of reflective $\omega$-sequences may have modifications that are applicable to some models of V=HOD.

**Theorem 15: (a)** There is a $\Sigma_2^V$ formula $\varphi$ such that there is no ordinal definable set $S$ with $\exists s \in S \, (\varphi(S) \Leftrightarrow \varphi(S \backslash \{s\}))$.
**(b)** There is a $\Sigma_2^V$ formula $\varphi$ such for every set of ordinals $S$ of order type $\omega$ with $\varphi(S) \Leftrightarrow \varphi(S \backslash \{S_0\}) \Leftrightarrow \varphi(S \backslash \{S_0, S_1\})$, every ordinal definable set includes (as a subset) or is disjoint from $S$ except for a finite set.
**Proof: (a)** Let $\varphi(S)$ be the following: $S$ is ordinal definable and $|T \backslash S \cup S \backslash T|$ is odd, where $T$ is the OD-least set such that $T \backslash S \cup S \backslash T$ is finite.

**(b)** Let $\varphi(S)$ be the following:
(1) There is an ordinal definable set $P$ such that both $S \cap P$, and $S \setminus P$ are infinite, and
(2) $P(S_0)$="$m$ is even", where $P$ is the HOD-least set that satisfies (1), and $m \geq 0$ is the minimum index such that $\neg(P(S_m) \Leftrightarrow P(S_{m+1}))$.
If (1) holds, then P is the same set for $S$, $S \setminus \{S_0\}$, and $S \setminus \{S_0, S_1\}$, and by applying (2) we can contradict $\varphi(S) \Leftrightarrow \varphi(S \setminus \{S_0\}) \Leftrightarrow \varphi(S \setminus \{S_0, S_1\})$. Thus (1) fails as required.

The conjecture (together with the proof of the theorem) implies that for a reflective sequence $S$, $\mathrm{HOD}(V_{\sup(S)})$ does not satisfy the axiom of choice, and every set in $\mathrm{HOD}(V_{\sup(S)})$ contains (as a subset modulo a finite set) or is disjoint from $S$ modulo a finite set.

Certain weak variations on the notion of $\omega$-sequences are well-understood (see the theorem below), and they appear to satisfy the conjecture.

Let $X$ be a predicate definable in $(V, \in, R_\alpha)$ from a set $S$ for some ordinal $\alpha$. A cardinal $\kappa$ is $X$-reflective iff the theory of $(V, \in, X, \kappa)$ with parameters in $V_\kappa$ is the same as the theory $(V, \in, X, \lambda)$ where $\lambda > \max(\kappa, \mathrm{rank}(S))$ and $R(\alpha+1, \lambda)$ holds. Iteration of reflectiveness are defined analogously. Note that $X$ can be coded by a set, and the notion of $X$-reflective (and $X$-$\beta$-reflective) is formalizable in our language.

The following theorem states that definability using $\omega^2 \alpha$ iterations of reflectiveness approximately corresponds to definability from a reflective $\omega$-sequence $S$ and $\alpha$ levels of the constructible hierarchy above $\{V_{\sup(S)}, S\}$.

**Theorem 16** (B1-B8): Let $S = \{S_i : i < \omega\}$ be an ascending sequence of $S_0$-reflective cardinals, and $X$ definable in $(V, \in, s, R_a)$ where $s$ is a set, $a$ is an ordinal, and $S_0 > \max(a, \mathrm{rank}(s))$. Let $\lambda$ be $\sup(S)$. Let $P$ denote the power set operation.
**(a)** The property of being an $X$-$\omega$-reflective cardinal is $\Pi_1^{P,X,S}$ definable in $V_\lambda$.
**(b)** Conversely, if $T$ is an ascending sequence of $X$-$\omega$-reflective cardinals with $T_0 \geq S_0$ and $\sup(T) \leq \lambda$, then $\Sigma_1^{P,X,S}$ theory of $V_\lambda$ with parameters in $V_{S_0}$ agrees with the corresponding $\Sigma_1^{P,X,T}$ theory of $V_{\sup(T)}$.
**(c)** (b) with $\Sigma_1$ definability replaced by $\Sigma_{<\omega}$ definability and elements of $T$ $X$-$\omega^2$-reflective.
**(d)** (b) with $\Sigma_1$ definability replaced by definability in $(L_\alpha(V_\lambda, S, X \cap V_\lambda), \in, S, X \cap V_\lambda)$ (and $(L_\alpha(V_{\sup(T)}, T, X \cap V_{\sup(T)}), \in, T, X \cap V_{\sup(T)}))$ where $\alpha < S_0$ and elements of $T$ are $X$-$\omega^2 \alpha + \omega^2$-reflective.
**(e)** The property of being an $X$-$\omega^2 \alpha$-reflective cardinal ($\alpha \leq S_0$) below $\lambda$ is definable in $(L_\alpha(V_\lambda, S, X \cap V_\lambda), \in, S, X \cap V_\lambda)$.

**Proof:** The proof is a routine induction on levels of reflectiveness, with (a) and (b) being the key inductive step.

We could continue iterations of reflectiveness beyond Ord by relying on an ordinal notation system. However, the proper course is to switch infinite sequences.

Let $M$ be a model of ZFC + "every set has a sharp" that has a measurable cardinal $\kappa$. Let $U$ be a normal ultrafilter on $\kappa$, $j_\alpha$ the corresponding $\alpha$-iterated ultrapower embedding, $M^\omega = j_\omega(M)$, and $S = \{j_\alpha(\kappa) : \alpha < \omega\}$. Note that $S$ is Prikry generic over $M^\omega$. Let $R$ — defined in $M$ — be the set of $x$ below $\kappa$ such that the theory of $(L(V^{M^\omega}_{j_\omega(\kappa)})[S], \in, x \cup S)$ with parameters in $V_x$ agrees with the theory of $(L(V^{M^\omega}_{j_\omega(\kappa)})[S], \in, S)$.

**Question:** Does $R$ — and the set of subsets of $R$ of order type $\omega$ — satisfy the basic properties for the weakening of reflective $\omega$-sequences to definability in the constructible universe above the supremum of the sequence?

If yes, and if our construction is nonrestrictive (that is its axiomatization does not prove false statements about the weakening of reflective $\omega$-sequences), then we could try to go further by using Prikry-Magidor forcing to get sequences longer than $\omega$.

### 7.3 Local Versions of Reflectiveness

The concept of reflectiveness can also be applied locally, below a regular uncountable cardinal. By the vastness of $V$, there should be enough space to have a notion of reflective ordinals there. Conceivably, using reflective sequences, one could have a local — that is in $V_\kappa$ for some small $\kappa$ — definition of the truth predicate for set theory, but the more likely possibility is that such sequences are definable not far beyond $V_\kappa$, and thus can be studied in the ordinary language of set theory.

For countable ordinals, the key result is that assuming sufficient determinacy, all sufficiently fast growing sequences of Turing degrees of length $\omega_1$ are effectively indistinguishable (Taranovsky 2012). A Turing degree $T$ naturally maps to an ordinal, specifically, the superior of $T$-recursive well-orderings. Thus, it appears that the predicate for lightface reflective for third order arithmetic sequences of countable ordinals of length $\omega_1$ that exclude their limit points is definable in third order arithmetic from the satisfaction relation for third order arithmetic.

In some cases, there is a set of incompatible properties, and to define the analogue of reflectiveness, we would have to choose one of the properties. For example, below the least inaccessible cardinal, to define an analogue of reflective cardinals, we would have to choose their cofinality. However, this diversity is not random. For every regular uncountable cardinal, there is a well-ordered set (modulo nonstationary sets) of canonical stationary sets. Once we pick a canonical ordinal definable stationary set, the analogue of reflectiveness appears to become unique — while some notions are stronger than others, all appear to be compatible. We prove a strong result in that direction using a plausible determinacy assumption (whose consistency, however, is not known).

**Definition:** *Extended Determinacy Maximum* is the assertion that ordinal definable (in the sense below) two player perfect information games of arbitrary ordinal length with positions being sets definable from countable sets of ordinals are determined.
Here, ordinal definable means that there is ordinal definable:
- initial position
- move function: position, player number, move → new position
- limit function $f$ (applied at all limit stages): $f(p)$ → new position, where $p$ is a sequence of positions corresponding to a valid play (starting at the initial position) such that for every valid $s$ and $t$ of the same length $f(s) = f(t)$ holds if for cofinally many $\alpha <$ length($s$), $s(\alpha) = t(\alpha)$.
- payoff set (set of positions won by the first player).

**Theorem 17** (Extended Determinacy Maximum): Let $\kappa$ be a regular uncountable cardinal such for every $\forall \lambda < \kappa \, \lambda^{\aleph_0} < \kappa$. Let $A \in \mathrm{OD}(\mathrm{Ord}^\omega)$ be a set of subsets of $\kappa$ of order type $\omega_1$, and $S \in \mathrm{HOD}(\mathrm{Ord}^\omega)$ be a stationary subset of $\kappa$. There is a set $X \subset S$, and a function $f: \kappa \to \kappa$ such that (1) for every $g: \kappa \to \kappa$ with $\forall x < \kappa \, g(x) \geq f(x)$, $g^\alpha(0) \in X$ where $\alpha$ is the least such that $g^\alpha(0) \in S$, and (2) $X$ is homogeneous for $A$, that is $A(x)$ is independent of $x$ for all $x \subset X$ of order type $\omega_1$. Furthermore, the above holds if $A \in \mathrm{OD}(\mathrm{Ord}^\omega)$ is a function with $A(x) \subset \min(x) \cup \omega$, with homogeneity defined by $A(x) = A(y) \cap A(x)$ where $\min(x) \leq \min(y)$.
**Proof:** Consider the following game: The players pick ordinals below $\kappa$, each of which must be above all previous ordinals, and let $p$ be the resulting sequence. Set to $x$ to be the set of limit points of $p$ below $\sup(p)$ that are in $S$. The game ends when $|x| = \omega_1$, and the first player wins iff $A(x)$. By stationarity of $S$, $\kappa$ will not be reached before $|x|$ will reach $\omega_1$. Given a winning strategy $T$, there is a strategy $T'$ that depends only on $x$ and the last ordinal. (Such $T'$ exists because the payoff depends only this information, and the last ordinal played is increasing. We construct such $T'$ from $T$ by recursively mapping positions consistent with $T$ into plays consistent with $T$, and doing this in a consistent way.)
By the closure of $\kappa$, there is a function $f$ such that if all ordinals played so far are $\leq a$, then the move by $T'$ is below $f(a)$. $f$ is the desired function. The rest of the theorem easily follows.
**Note:** Determinacy for games of length $\omega_1$ suffices for S having only elements of cofinality $\omega$.

An interesting question is whether under appropriate assumptions, the theorem extends to $\{A, S\} \in \mathrm{OD}(\mathrm{Ord}^{<\kappa})$ and sequences of length $\kappa$. A consequence of the theorem is that the analogues of reflectiveness cannot be used to define real numbers that are not ordinal definable. True reflective sequences can also be defined below a cardinal $\kappa$ that has sufficiently strong reflective properties. The difference is that the relevant stationary set is not ordinal definable.

### 7.4 Further Extensions

**Note:** This is the original version. A new exposition with more details is in

"Reflective Sequences" below.

Continuing our quest to higher levels of expressiveness beyond reflective $\omega$-sequences:

**Notion 2:** For any ordinal $\alpha$, a reflective sequence $S$ of length $\alpha$ with $\min(S) > \alpha$. $S$ excludes its limit points.

**Notion 3:** A reflective sequence $S$ with the order type of $S$ having sufficiently strong reflection properties relative to $S$. $S$ excludes its limit points.

**Notion 4:** A reflective sequence $S$ with length extremely large relative to $\min(S)$. $S$ excludes its limit points except for points of extremely high reflectiveness. $S$ has a maximum element, which, in terms of limit points, has maximal reflectiveness.

In Notion 4, a question arises about which limit points to include. We clarify the type of $S$ as follows: There is a transitive model $M$ of ZFC such that $S$ is the set of measurable cardinals in $M$ and $\max(S)$ is the only measurable cardinal $\kappa$ of Mitchell order 2 (that is, $\kappa$ has a normal measure that concentrates on measurable cardinals).

For each of these notions, there is a corresponding property of ordinals: $R_{\text{ord}}(\kappa)$ iff $\kappa = \min(S)$ for some $S \in R$. For Notions 1-3, it might be possible to recover the notion by using sequences of the appropriate type of ordinals in $R_{\text{ord}}$; for Notion 4, we would also need $R_{\text{ord}}$-reflective cardinals (or a sufficiently expressive weakening) for recovery of the notion from $R_{\text{ord}}$. An interesting conjecture is that if we weaken Notion 4 by using $(V_{\sup(S)}, \in, S)$ in place of $(V, \in, S)$ for testing correctness of $S$, then the resulting $R_{\text{ord}}$ is such that $R_{\text{ord}}(\kappa)$ holds iff the theory of $(V, \in, (R_{\text{ord}} \backslash \kappa) \cup \{\kappa\})$ with parameters in $V_\kappa$ agrees with the theory of $(V, \in, R_{\text{ord}} \backslash \lambda)$ for every $\lambda \geq \kappa$.

We can go beyond Notion 4 by including more measurable cardinals in $M$, and we conclude with one of the most expressive extensions to the language of set theory that is currently known.

**Notion 5:** Define $R(f)$ to hold iff
(1) There is a transitive model $M$ of ZFC such that $\text{domain}(f)$ is the set of measurable cardinals in $M$, and for every $\kappa$ in the domain of $f$, $f(\kappa)$ is the Mitchell order of $\kappa$ in $M$. $M$ has only one cardinal $\kappa$ such that $o(\kappa)^M = \kappa + 1$, and it is the largest measurable cardinal in $M$.
(2) For the class functions described above, $f$ satisfies reflecting and large cardinal properties so strong that the theory of $(V, \in, f)$ with parameters in $V_{\min(\text{domain}(f))}$ is independent of $f$.
[Note: We are not sure how to best handle parameters between $\min(\text{domain}(f))$ and $\max(\text{domain}(f))$, but our handling probably does not diminish the expressive power of the predicate.]

We can try to go further by making $M$ an iterate of the minimal inner model with a strong cardinal (and $\text{domain}(f)$ the set of measurable cardinals in $M$ below the strong cardinal), but I am not sure that this notion is well-defined, so for now, we will end our adventure here. As the expressiveness is increased, the notions become more technical, and the question is: Is there a simpler and fundamentally more expressive extension to the language of set theory?

# 8 Reflective Sequences

**Note:** This section was written as a standalone paper, but per reader's convenience and request is also incorporated here. This section can be read independently of the previous sections.

## 8.1 Introduction

Within the language of set theory, one reaches higher and higher expressive power by climbing higher in the cumulative hierarchy of V. But how can we go further once the language allows quantification over the whole V? Intuitively, we would want to continue the hierarchy above V, except that all sets are already in V. The solution is to label a cardinal κ such that $V_κ$ is sufficiently close to V, and continue the hierarchy above $V_κ$. $V_κ$ represents V, and with κ labeled in an extended language, hierarchy above $V_κ$ corresponds to higher order set theory.

To go further, we can iterate higher order set theory by picking λ>κ with $V_λ$ representing V, and then μ>λ, and continuing to longer sequences of ordinals. Thus, continuing the cumulative hierarchy above V corresponds to labeling certain ordinals that are sufficiently similar to Ord — in other words, certain ordinals with sufficiently strong reflection properties — while staying within V.

How do we choose the right κ? The answer is that we postulate a certain degree of symmetry and reflection in V.
**Convergence Hypothesis (general form):** For an appropriate type of objects, all objects of that type with sufficient reflection properties are, to a certain degree, indistinguishable from each other.
**Convergence Hypothesis (ordinals):** If α and β are ordinals with sufficiently strong reflection properties, then φ(S, α) ⇔ φ(S, β) whenever rank(S) < min(α, β) and φ is a first order formula of set theory with two free variables.
**Definition:** κ is a *reflective cardinal*, denoted by R(κ), iff (V,∈,κ) has the same theory with parameters in $V_κ$ as (V,∈,λ) for every cardinal λ>κ with sufficiently strong reflection properties.
**Example:** In L (assuming $0^{\#}$), every Silver indiscernible has sufficient reflection properties, which allows us to define $R^L$ in V.
**Note:** A formalist can treat the Convergence Hypothesis as a guiding principle to get good systems of axioms. A formalist can ask, for example, whether we want the Convergence Hypothesis to hold in V, and if yes (or maybe), investigate how far it can hold and how to axiomatize it. A counterexample to the Convergence Hypothesis would be a property that 'toggles', with an explanation of why it toggles.

To the extent that it holds, the Convergence Hypothesis allows us to define new reflective notions without ambiguity. We specify a notion by stating its type (example: a pair of ordinals) along with the class of predicates for which it agrees with all objects (of the same type) with sufficient reflection properties.

**Definition:** Given a length condition P (such as having length ω), an increasing sequence of ordinals S satisfying P is *reflective* (denoted by $R_P(S)$ or simply $R(S)$) iff
   (1) the theory of $(V,\in,S)$ with parameters in $V_{\min(S)}$ is correct, that is it agrees with $(V,\in,T)$ for every T satisfying P and having sufficient reflection properties and $\min(T) > \min(S)$, and
   (2) for every $\alpha < \sup(S)$ $S\backslash\alpha$ is reflective under the above criterion (where the length condition is modified to correspond to $S\backslash\alpha$ if $S\backslash\alpha$ is "shorter" than S).
A variation on the notion is to omit the second condition.

*Notes about R:*
• We use R rather a single element of R for two reasons:
  - We want to be able to express predicates, such as the satisfaction relation of $(V,\in)$, and not just statements.
  - For infinite reflective sequences, it is unclear if (or how) an example is definable, and we want the intended semantics of the language to be unambiguous.
• There is a hierarchy notions related to R based on the degree of "correctness" required (one example is to allow just $\Sigma^V_2 \varphi$ in the definition as opposed to $\Sigma^V_{<\omega} \varphi$). Essentially, our choice for R(x) is to require just enough correctness to maximize expressiveness of $(V,\in,x)$ for predicates that do not depend on x.
• In every extension here, R will correspond to a single type of objects (for example, sequences of ordinals of length ω). However, these notions form a hierarchy, and an alternative, which we will use informally, is to use R(x) in place of $R_{\text{type}(x)}(x)$, thus having just a single symbol R.

## 8.2 Reflective Sequences

To understand the hierarchy of reflective sequences, imagine a process of picking ordinals with sufficient reflection properties, and with the process itself having sufficient reflection properties. Let S be a resulting sequence. Here are some stages in the process:
**1.** n ordinals for finite n. This is well-understood (see [Iterations of Reflectiveness](#) above for theory).
**2.** ω ordinals. See [Reflective ω-Sequences](#) above for some results. There is a refined hierarchy based on the degree of correctness we require. We have good axiomatization (specifically, through iterations of reflective cardinals) up to about correctness in $L_{\sup(S)}(V_{\sup(S)})$, and through an analogue of Prikry forcing, we may be able get up to $M(V_{\sup(S)})$ (where M is a mouse operator and assuming sufficient structure above M). However, under reasonable assumptions, full correctness in V contradicts V=HOD, which makes analysis and finding of canonical models difficult. A reasonable conjecture is that the $\Sigma_2$ truth predicate of $(V_{\sup(S)}, \in, S)$ corresponds to enumeration of all ordinal definable reals.
**3.** α ordinals for countable α. We expect R to be closed under subsequences that preserve order type.
**4.** $\omega_1$ ordinals. Assuming the Continuum Hypothesis (CH), we cannot require closure of R under all subsequences of the same length. Even without CH, this limitation holds at $\omega_2$.
**5.** α ordinals for a fixed ordinal $\alpha > \omega_1$ ($\min(S)>\alpha$).
**6.** More complicated conditions about the length of S, such as $|S|=\min(S)$. S excludes its limit points.
**7.** Length of S has sufficiently strong reflection properties relative to S. S excludes its limit points.
**8.** To go beyond (7), we assume that S includes its limit points at places with sufficient (relative to S) degree of reflectiveness. One may then state conditions on the length of S such as "S includes exactly ω limit points, and they are cofinal in S" or "S has maximum element and it is the only element such that S below it is stationary".
**9.** Length of S has sufficiently strong reflection properties relative to S. This notion is axiomatized in "Axiomatization of Long Reflective Sequences" below.

To go beyond what is expressible with S, we mark some points on S that have sufficient (relative to S) degree of reflectiveness, analogously to how S is obtained by marking some ordinals with sufficient reflectiveness in V. To go further, we then add a third type of marking, and in general, for each ordinal in S, we can assign a degree of reflectiveness through a function $f:S \to Ord \setminus \{\emptyset\}$. For a limit α, $f(x)=\alpha$ corresponds to x receiving all markings <α. In the definition of reflectiveness, use f in place of S, min(dom(f)) in place of min(S) (and min(dom(T)) in place of min(T)), sup(dom(f)) in place of sup(S), and in place of S\α, use $f \restriction (Ord \setminus \alpha)$ modified to make f(min(dom(f)\α)) equal 1. Here are some notions. (Here S=dom(f).)

**10.** For a fixed ordinal α, $f(sup(S))=\alpha$ (and hence sup(S)∈S), and sup(S) is the only ordinal κ with $f(\kappa)=\alpha$ (and min(dom(f))>α).
**11.** $f(sup(S))=sup(S)$, and sup(S) is the only ordinal κ with $f(\kappa)=\kappa$.
**12.** To go beyond $f(\kappa)=\kappa$, we modify the condition — that $f(\kappa)=\alpha+1$ implies that κ has sufficiently strong reflection properties relative to f that is clipped to α — by invoking an ordinal notation system (for ordinals ≥κ) or just a comparison method between f(λ) and f(κ), and considering f clipped below f(κ) as defined through the comparison method. Below (in "Extension to Reflective Functions") we axiomatize the following: There is only one ordinal κ with $f(\kappa)=(\kappa^+)^{HOD}$, and it equals sup(dom(f)).

For reflective cardinals, an axiomatization essentially states the degree of correctness (relative to other reflective cardinals) and existence of enough reflective cardinals, and it can be extended with large cardinal axioms applied to reflective cardinals. For reflective sequences, we also need to state enough axioms for the sequence as a whole, or the axiomatization would be too incomplete. Just like large cardinal axioms for cardinals, there are various axioms for reflective sequences.

The type of f and many of its basic properties can be formalized using inner models. For example, for notion 11 we can require that there is an inner model M of ZFC such that
 - dom(f) is the set of measurable cardinals in M
 - for every κ in the domain of f, f(κ) is the Mitchell order of κ in M
 - M has only one cardinal κ such that $o(\kappa)^M=\kappa$ and it is the largest measurable cardinal in M.
From inner model theory, it appears that without restricting the class of f that satisfy the requirement, M can be taken to be an iterate of the minimal inner model with these cardinals. The possibilities for M and iteration are sufficiently rich to get any reasonable f corresponding to this notion (notion 11 in this case). One can also similarly state requirements on M for other notions.

If f satisfies the condition above (with M an iterate of the minimal model), then the theory of (L(f),∈,f) is independent of f, which is one example of the Convergence Hypothesis, and moreover, we can study this theory to get a better understanding of f.

Our axiomatization involves a strengthening of the above condition by asserting that elements of a reflective sequence are, in a sense, infinitary indiscernibles. To do that, we need to know how much infinitary indiscernibility is possible in V.

## 8.3 Limits on Infinitary Indiscernibility

The axiom of choice sharply limits the amount of infinitary indiscernibility.

**Theorem 18:** For every uncountable sequence of distinct reals X, there is X-definable f:V→ω such that for every S⊂Ord of order type $\omega_1$ there is X,S-definable T⊂S of order type $\omega_1$ with f(T)≠f(S).
**Notes:**
• This implies that for every κ≥$\omega_1$, the partition relation κ → $\omega_1^{\omega_1}_\omega$ fails for a function f definable from X.
• For definability, $\Delta^V_2$ definability more than suffices.
• If T could be made S-definable, the axioms below would be inconsistent, but a weakening to (essentially) lightface version would be unaffected.
• *Acknowledgment:* I would like to thank Hugh Woodin for outlining a proof of a weaker version of the theorem with exponent $\omega_2$.
**Proof:**
Using X, we can define a sequence X′ of subsets of ω of length $\omega_1$ that does not have a perfect subset. Thus, the following game (this is the standard game for the perfect subset property) is undetermined:
Player I: Each move consists of zero or more zeroes and ones.
Player II: Each move is 0 or 1.
Player I wins iff if the concatenation of the moves belongs to the set (X′).
However, analogously to the standard proof of analytic determinacy from sharps (which can be found in, for example, Wikipedia [Determinacy](#) article (accessed October 30, 2016)), the partition relation would make the game determined as follows. The tree $T_{X′}$ corresponding to the game will have $\omega_1$ branches at the root, each corresponding to an element of X′, and no further branches. The open auxiliary game requires player II to play the subtree of $T_{X′}$ corresponding to the play so far (including all initial segments of the play), and a mapping of the subtree into ordinals, such that the mappings are consistent with each other and agree with the Kleene-Brouwer order. Player I wins the auxiliary game and has a definable strategy: For example, pick the least move that gives maximum reduction of the rank of the remaining game. Let g(S) return move given position where S is the set of ordinals used in the tree. If we attempt to use S to convert the strategy into a strategy for the original game, the failure will give us the required T. While g returns elements of $2^\omega$, it can be converted into f:V→ω as follows: f(S)=2n+(g(S))$_n$ where n = min(m: ∃T⊂S T∈Def(X′,S) ∧ |T|=$\omega_1$ ∧ (g(S))$_m$≠(g(T))$_m$)) where Def(X′,S) can be for example definability from X′,S in $V_{sup(S)+1}$.

**Corollary 19:** There is a definable predicate P such that for every set of ordinals S

of order type at least $(2^{2^{\omega_1}})^+$ there are subsets $T_1$ and $T_2$ of order type $\omega_1*2$ such that such that $P(T_1) \neq P(T_2)$).
**Proof Sketch:** Here is one such P: P(S) iff $\{S_i : i < \omega_1\}$ agrees with $\{S_i : \omega_1 \leq i < \omega_1*2\}$ with respect to all predicates that are definable in $V_{\sup(S)}$ using parameters that are subsets of $\omega_1$.

In the absence of a convincing theory, lack of closure under subsets would cast doubt on whether the notions of reflective sequences are well-defined. However, we can weaken closure under subsets by requiring the subsets to be definable, and that may give us what we need to axiomatize R.
**Note:** It is also possible that the convergence hypothesis holds for statements but not for statements with arbitrary sets as parameters. In that case, R would be ill-defined in V but might have a preferred well-defined theory, and its lightface version (requiring correctness only for statements without parameters) would be well defined. If the 'boldface' version is ill-defined or otherwise undesirable, the systems of axioms below can be modified for the lightface version by requiring the relevant theories to agree without parameters (as opposed to with parameters), and analogously for a strengthening of the lightface version that allows countable sequences of ordinals as parameters.

On the other hand, assuming determinacy of games on ordinals of length $\omega_1$ and ordinal definable payoff, $\kappa \to \kappa^{\omega_1}_{2^\lambda}$ holds for every regular countably closed $\kappa$ and $\lambda < \kappa$ for functions definable from a countable sequence of ordinals. See [Local Versions of Reflectiveness](#) above for a proof sketch of a more general proposition using Extended Determinacy Maximum. *Question:* Can similar relations can be obtained from much weaker assumptions?

Partition relations with exponents $>\omega_1$ might be consistent by restricting to definable subsets. A strengthening of the strong partition relation (restricted to definable partitions and subsets) is to also require the homogeneous set to have enough limit points and require the same of the subset, and this can be strengthened further by requiring the sets to be in an appropriate filter. Here is a key proposition to investigate:
**Conjecture:** There is an ordinal $\kappa$, a normal nonprincipal $\kappa$-complete ultrafilter U on $\kappa$, and $S \in U$ such that for every $T \subset S$ with $T \in U$ and T ordinal definable from S and a countable set of ordinals, $(V, \in, T)$ and $(V, \in, S)$ have the same $\Sigma^V_2$ theory with parameters in $V_{\min(S)}$.
*Question:* Is this conjecture consistent?
If consistent and nonrestrictive, this property can be used in axiomatization of long reflective sequences of ordinals as follows.

## 8.4 Axiomatization of Long Reflective Sequences

Here is an axiomatization corresponding to notion 9 above. Its consistency and truthfulness are unclear.

**Language:** $\in$, R (unary predicate)
*Definition:* Theory(M, $\in$, S; W) (where M is transitive, $S \subset M$ and $W \subset M$) is $\{$"$\varphi$",w: $w \in W \land (M, \in, S) \vDash \varphi(w)\}$ where $\varphi$ has 1 free variable and may use relational symbols '$\in$' and 'S' (unary predicate).

*Definition:* $OD_S(C)$ is the class of sets ordinal definable (without using R) from S and elements of C.
*Definition:* $R_1(\kappa) \Leftrightarrow \exists S \in R\ \kappa \in S$
*Note:* sup stands for supremum.
**Axioms:**
**1.** ZFC
**2.** $\forall x \exists y \forall z (z \in y \Leftrightarrow z \in x \wedge R(z))$
**3.** $R(S) \Leftrightarrow R_{aux}(S) \wedge \forall \alpha < \sup(S)\ R_{aux}(S \backslash \alpha)$

where $R_{aux}(S) \Leftrightarrow S \neq \emptyset \wedge S \subset Ord \wedge \forall T \in R\ \forall \kappa \in R_1 \backslash (\sup(S \cup T)+1)\ Theory(V_\kappa, \in, S; V_{\min(S \cup T)}) = Theory(V_\kappa, \in, T; V_{\min(S \cup T)})$
**4.** (schema; n is a natural number) $\exists \kappa \in R_1\ \exists S \in R\ \kappa = \sup(S) \wedge V_\kappa <_{\Sigma^R_n} V$
**5.** $\exists S \in R\ \exists U \ni S$ (a normal nonprincipal κ-complete ultrafilter on $\kappa = \sup(S)$)
$\forall T \in OD_S(Ord^\omega)\ (T \subset S \wedge T \in U \Rightarrow R(T))$.

*Notes:*
- (3) specifies the type of R and the degree of correctness required. The intended semantics is that members of R have sufficiently strong reflection properties so as to (if possible) render R unambiguous in V. As mentioned above, the general underlying hypothesis is that for an appropriate type, all objects of that type with sufficiently strong reflection properties are sufficiently indistinguishable from each other.
- $R_{aux}$ and R have the same expressive power (the definition of $R_{aux}$ in (3) could just as well have used $R_{aux}$; $R_1$ is $\{\min(S): S \in R_{aux}\}$), but R is more natural (and is different from $R_{aux}$). Other than that, (1)-(4) are essentially standard, and the combinatorial properties are obtained from (5).
- (4) is a strengthening of replacement for formulas with R, and it also ensures that sup(S) has essentially the same properties as members of S.
- A weakening of (5) (if one objects to requiring measurable cardinals in V) is to require U to be an ultrafilter on a smaller model ZFC, as long as it contains $HOD(V_\kappa, S)$.
- It is unclear whether the property in (5) holds for all members of R. Its weakening to $\Sigma^V_n$ predicates holds, but it is unclear whether U can be chosen independent of n.
- In (5), it is unclear whether the ultrafilter itself can be chosen to have a canonical theory with parameters in $V_\kappa$.
- A strengthening of 5 is to also require $S \in i^U(R)$ to hold (where $i^U$ is the embedding corresponding to U).
- These notes also apply to the extension to reflective functions below (except the part about $S \in i^U(R)$).

Using the theorem above, the axioms imply that there is no $OD(Ord^\omega)$ uncountable sequence of distinct reals.

*Question:* Is it consistent/provable/true that $\exists \kappa > 0\ \forall f: \kappa \to \kappa\ \exists T \in R \cap P(\kappa)\ \forall \alpha < \kappa\ T_{\alpha+1} > f(T_\alpha)$?
If yes, what about the stronger $\forall S \in U\ \exists T \subset S\ T \in R \cap U$ (where U is as in axiom 5 above)?
*Notes:*
* The same question can be asked about extensions to reflective functions using the

domain of the function in place of T.
\* The statement (either weaker or stronger form) is natural and appears to fill in the gap created by the lack of closure of R under subsequences. However, we do not know whether there are 'pathological' counterexamples to the statement.
\* If true but unprovable, the statement is a good addition to the axioms.

Here are possible outcomes for the axioms:
a. The axioms are inconsistent. One would then search for a suitable replacement, and if it is found, try again.
b. The axioms are consistent but restrictive (for example, they contradict a true large cardinal axiom). One would then also look for replacement (though models for the axioms may also be interesting).
c. The axioms are consistent but too incomplete because the underlying notion is ill-defined. One can still study models with the axioms, and appreciate the underlying symmetry and beauty.
d. The axioms are consistent but too incomplete. One would then search for additional axioms. Good axiomatization is key evidence that the underlying notion is well-defined.
e. The axioms form the core of the theory, but it can still be extended, analogously to ZFC forming basic axioms for set theory.

## 8.5 Extension to Reflective Functions

We can go beyond notion 9 by having a coherent system of measures instead of a single measure. For an ordinal $\kappa$, measure $U(\kappa,\alpha)$ (note that its Mitchell rank need not be $\alpha$) will witness that there are enough ordinals $\lambda<\kappa$ with $f(\lambda)=\alpha$ (or if $\alpha \geq \kappa$ that $f(\lambda)$ corresponds to $\alpha$). For transformations of f, besides reducing its domain, we can downgrade some of the labels. A complication is that in the transformed f, we may have, for example, an ordinal with label 2 without enough label 1 ordinals below it. We resolve it by recursively relabeling ordinals, each time keeping as much structure as is compatible with the labels for lower ordinals. The proposed axiomatization for notion 12 is as follows.

**Language:** $\in$, R (unary predicate)
*Definition:* $\kappa \in R_1 \Leftrightarrow \exists f \in R$ "f is a nonempty function with $\text{dom}(f) \subset \text{Ord}$" $\wedge$ $\kappa = \min(\text{dom}(f))$.
*Note:* See above for definitions of Theory(M, $\in$, S; W) and OD.
**Axioms:**
**1.** ZFC
**2.** $\forall x \exists y \forall z (z \in y \Leftrightarrow z \in x \wedge R(z))$
**3.** $R(f) \Leftrightarrow R_{aux}(f) \wedge \forall \alpha < \sup(\text{dom}(f)) \exists \beta > \alpha\ R_{aux}(f \restriction (\text{Ord} \backslash \beta))$
   where $R_{aux}(f) \Leftrightarrow f \neq \varnothing \wedge$ (f is a function) $\wedge \text{dom}(f) \subset \text{Ord} \wedge \text{ran}(f) \subset \text{Ord} \wedge \forall g \in R$ $\forall \kappa \in R_1 \backslash (\sup(\text{dom}(f) \cup \text{dom}(g))+1)$ Theory($V_\kappa$, $\in$, f; $V_{\min(\text{dom}(f) \cup \text{dom}(g))}$) = Theory($V_\kappa$, $\in$, g; $V_{\min(\text{dom}(f) \cup \text{dom}(g))}$).
**4.** (length condition) $\forall f \in R\ \{\kappa \in \text{dom}(f): f(\kappa) = (\kappa^+)^{HOD}\} = \{\sup(\text{dom}(f))\}$
**5.** (schema; n is a natural number) $\exists \kappa \in R_1\ \exists f \in R\ \kappa = \sup(\text{dom}(f)) \wedge V_\kappa \prec_{\Sigma^R_n} V$
**6.** (key combinatorial principle) There is $f \in R$ and a sequence of ultrafilters U such that:
   i. $\text{dom}(U) = \{(\kappa,\alpha): \kappa \in \text{dom}(f) \wedge \alpha < f(\kappa)\}$
   ii. $U(\kappa,\alpha)$ is a $\kappa$-complete normal nonprincipal ultrafilter on $\kappa$

   iii. ∀(κ,α)∈dom(U) (κ,α)∉dom($i^{U(κ,α)}$(U)) ∧ ∀β<α $i^{U(κ,α)}$(U)(κ,β) = U(κ,β)
(coherence condition; *i* is the embedding corresponding to the ultrafilter)
   iv. Given a function g:Ord→Ord\{∅} with dom(g)⊂dom(f), define g′ (with dom(g′)=dom(g)) as follows: g′(κ)=min(g(κ), sup(α+1: α=0 ∨ ∃β>0 $i^{U(κ,β)}$(g′)(κ) = α)). *Key Condition:* ∀g∈$OD_f$($Ord^ω$) g′(sup(dom(f))) = f(sup(dom(f))) ⇒ R(g′).

## 8.6 Further Extensions

### General

To extend the language further, some of the possibilities are:
1. Extend the length condition to a higher ordinal.
2. Extend R to a different type of objects.
3. Find a different property in the set-theoretical universe.
4. Find something other than sets that is not reducible to sets.

The axioms above are essentially independent of the length condition, and by using (non-overlapping) extenders in place of ultrafilters, we can go up to a strong cardinal (though the resulting system might well be inconsistent or otherwise false). On the other hand, our semantic development depends on being able to compare f(κ) with f(λ). We can likely go beyond $(κ^+)^{HOD}$ by, for example, using R to define a higher ordinal and requiring correctness relative to R in the axioms to define R'. To potentially go much further, note that we depend only on the correctness of theory with appropriate parameters rather being able to compare all f(κ) and f(λ), and it is possible that there is unique privileged theory here.

*Question:* Does the convergence hypothesis apply to the class of iterates of the minimal inner model with a strong cardinal?

If yes, then we have a corresponding extension to the language of set theory, essentially corresponding to maximal canonical non-overlapping f.

### Overlapping Extenders

We can try to go even further by allowing f to overlap, that is by letting f(κ)≥λ for λ>κ and λ∈dom(f). This may be analogous to overlapping extenders in the models beyond a strong cardinal — for example two cardinals strong up to the same measurable cardinal may correspond with f(κ) = f(λ) = min(dom(f)\(λ+1)) (with κ<λ) — but there are other possibilities as well, and the right framework is unclear.

For overlapping f, we must impose a restriction on the type of overlap. Otherwise, let $x_0$<$x_1$<$x_2$<… be an infinite sequence of ordinals such that ∀i<ω f($x_i$) ≥ sup($x_i$:i<ω). Let $y_i$=min(α: f(α)≥$x_i$). Without restrictions on overlapping, we would have $y_0$<$y_1$<…<$x_0$, and repeating the process with y in place of x, we get $z_0$<$z_1$<…<$y_0$, and so on, leading to infinite regress. One restriction is to prohibit infinite overlap — that is to prohibit $x_0$<$x_1$<… with ∀i<ω f($x_i$)≥$x_{i+1}$. However, this restriction is too limiting, and there is a better alternative. Note that if κ < λ < μ and κ is <λ-strong and λ is <μ strong, then κ is <μ strong. Analogously, we require that ∀κ∈dom(f)∀λ∈dom(f) (κ<λ ∧ f(κ)≥λ ⇒ f(κ)≥f(λ)).

The right general framework is to have a rich system of embeddings to witness f, perhaps using iteration trees. For an example, suppose $f(\kappa)=\lambda+1$ where $\lambda=\min(\text{dom}(f)\setminus(\kappa+1))$ and $\mu>\lambda$ and $\mu\in\text{dom}(f)$. If we remove all ordinals strictly between $\kappa$ and $\mu$ from dom(f), then we should be able to arrange $f(\kappa)=\mu+1$, and the system should have an embedding to witness that. However, we do not yet know how to axiomatize such systems. Below, we attempt an axiomatization for limited degree of overlap and without the ability to move $f(\kappa)$ upward. The later restriction might make the system too incomplete, but (if consistent) it should still show important structure.

**Language:** $\in$, R (unary predicate)
*Note*: See "Axiomatization of Long Reflective Sequences" above for definitions of Theory(M, $\in$, S; W) and OD.
**Axioms:**
**1**,**2**,**3** as above ("Extension to Reflective Functions").
**4.** (schema; $\underline{n}$ is a natural number) There is $f \in R$, $\kappa_{\max}=\sup(\text{dom}(f))$, a sequence of extenders E, and an ultrafilter U such that:
   i. $\text{dom}(E) = \{(\kappa,\alpha): \kappa\in\text{dom}(f) \wedge \alpha<f(\kappa)\}$
   ii. $E(\kappa,\alpha)$ is a normal short extender such that the corresponding embedding $i^{E(\kappa,\alpha)}$ is $\alpha$-dom(E)-S strong and has critical point $\kappa$. Here S is the satisfaction relation for $\Sigma^{R,V}_{\underline{n}}$ formulas.
   iii. $\forall(\kappa,\alpha)\in\text{dom}(E)\ (\kappa,\alpha)\notin\text{dom}(i^{E(\kappa,\alpha)}(E)) \wedge \forall\beta<\alpha\ i^{E(\kappa,\alpha)}(E)(\kappa,\beta) = E(\kappa,\beta)$ (coherence condition)
   iv. U is a $\kappa_{\max}$-complete normal nonprincipal ultrafilter on $\kappa_{\max}$.
   v. $V_{\kappa_{\max}} \prec_{\Sigma^R_{\underline{n}}} V$ and $\text{Theory}(V,\in,R,f;V_{\kappa_{\max}}) = \text{Theory}(i^U(V),\in,i^U(R),f;V_{\kappa_{\max}})$, where Theory is restricted to $\varphi$ with $\underline{n}$ unbounded quantifiers.
   vi. Given a function $g:\kappa_{\max}\to\text{Ord}$, define $g':\kappa_{\max}\to\text{Ord}$ as follows: $g'(\kappa) = \min(\{g(\kappa),$ $\sup(\alpha+1: \exists\beta\ i^{E(\kappa,\beta)}(g')(\kappa) = \alpha)\} \cup \{\sup(\alpha+1: \alpha<g'(\lambda) \wedge \forall\gamma<g'(\lambda)\exists\beta\ i^{E(\lambda,\beta)}(g')(\lambda)\geq\gamma \wedge i^{E(\lambda,\beta)}(g')(\kappa)\geq\alpha): \lambda<\kappa \wedge g'(\lambda)\geq\kappa\})$. Let $g'' = g'\restriction\{\kappa:g'(\kappa)>0\}$. *Key Condition:* $\forall g\in\text{OD}_{f,R}(\text{Ord}^\omega)\ (\{\lambda\in\text{dom}(g''): g''(\lambda)=\kappa_{\max}\}\in U \Rightarrow \text{Theory}(V,\in,R,g'';V_{\kappa_{\max}}) = \text{Theory}(V,\in,R,f;V_{\kappa_{\max}}))$, where Theory is restricted to $\varphi$ with $\underline{n}$ unbounded quantifiers. Furthermore, if $g\restriction\text{dom}(f) = f$, then $g''$ is f.

*Notes:*
- If the system is inconsistent, it may still be useful as an outline of how the axioms can look like.
- A schema was used to get the most natural stopping point at the price of some syntactic complexity. However, expressiveness using R goes only slightly beyond what is expressible with a single $f\in R$ (assuming the relevant objects are in $V_{\min(\text{dom}(f))}$), and there is a natural mild weakening of 4 that avoids the schema and goes just beyond what is expressible about f without using R:
   4.ii - let S be the satisfaction relation for (V,$\in$,f) formulas (which is definable using R and f). It is unclear whether omitting S entirely would cause much harm.
   4.v - simply require that $i^U(R)(f)$ holds, and as a partial substitute for the remainder, add axiom 5 (from "Extension to Reflective Functions") to the axioms (or a related statement/schema).
   4.vi - instead of using Theory, just require that $g''\in R$; also, prohibit use of R in defining g.
- $\alpha$-dom(E) strongness ensures that the embeddings witness not only the absolute

- largeness of f(κ) but also the structure of f between κ and f(κ).
- In 4.vi, note that if κ∉dom(f), then there is no E(κ,β), and sup(∅)=∅.
- In 4.vi, "g″ is f" is used to avoid having to write out some properties twice (once for f and once for g″).
- The difference of 4.vi from Extension to Reflective Functions 6.v is the treatment of overlapping f. If f(λ)>λ, then the structure of f between λ and f(λ) must be witnessed using E(λ,α) and hence by f below λ, and we apply the same requirement to g′. Thus, if the necessary structure below λ is erased, we correspondingly reduce g′ in the region between λ and f(λ).

For a model M with overlapping extenders, there may be multiple reasonable ways to represent M by f, and there is additional structure and complexity once M has a cardinal κ with a measure that concentrates on <κ-strong cardinals. One option is to let R apply to models directly (and for now, dropping the analogue of the second condition in the definition of R). Here is an example: R = {code(M): $M_0$ elementarily embeds into M and (V,∈,code(M)) has the correct theory with parameters in $V_κ$ where κ is the least measurable cardinal in M}, where $M_0$ is (for example) the minimal inner model with a measurable limit λ of <λ-S-strong cardinals where S is the set of <λ strong cardinals, and code(M)=$V_δ^M$ for the least δ such that M=L($V_δ^M$).

The study of reflective constructs has only just begun.